\documentclass[12pt]{amsart}

\theoremstyle{plain}
\newtheorem{theorem}{Theorem}[section]
\newtheorem{corollary}[theorem]{Corollary}
\newtheorem{lemma}[theorem]{Lemma}
\newtheorem{proposition}[theorem]{Proposition}
\newtheorem{example}[theorem]{Example}

\theoremstyle{definition}
\newtheorem{definition}[theorem]{Definition}

\theoremstyle{remark}
\def\qed{\hskip .6em \raise1.8pt\hbox{\vrule height4pt
width6pt depth2pt}}
\def\qedd{\hskip .4em \raise1.8pt\hbox{\vrule height3pt
width5pt depth1.8pt}}
\def\sq{\hskip .6em \raise1.8pt\hbox{\vrule
height4pt width6pt depth2pt}}

\def\R{\rm I\kern-2ptR}
\def\N{\rm I\kern-2ptN}

\font\titlos=cmb10 at 18 pt
\font\titloss=cmb10 at 16 pt

\begin{document}

\
\vskip -1.5cm

\centerline {\bf To Professor R. Cristescu on the occasion of his 70-th
birthday}

\vskip 1cm

%\hfill Draft, \today

\centerline {\titlos  d-independence and d-bases in vector lattices}

\bigskip

\centerline {\titloss Y. Abramovich and A. Kitover}

\bigskip
\section{Introduction}

This article contains the results of two types. 
In Section~3 we give a complete characterization of band preserving    
projection operators on Dedekind complete vector lattices. These
operators were instrumental in our work [AK2], and now we have obtained          
their description. This is done in Theorem~\ref{t:main.thm}.
Let us mention also  Theorem~\ref{t:char.pr.oper} that contains a    
description of such operators on arbitrary laterally complete vector   
lattices.
The central role in these descriptions is played by d-bases, 
one of two principal tools utilized in [AK2]. 
The concept of a d-basis, originally considered in this 
context in [AVK],
has been  applied so far   
only to vector lattices with a large amount of projection bands.             
The absence of the projection bands has been the major obstacle              
for extending, otherwise very useful concept of d-bases, 
to arbitrary vector
lattices. In Section~4 we will be able to overcome this obstacle by
finding a new way to introduce d-independence in an
arbitrary vector lattice. This allows us to produce a new
definition of a d-basis which is 
free of  the existence
of projection bands. We illustrate this by proving several
results  devoted to  cardinality of d-bases. 
Theorems~\ref{p:pr5} and~\ref{p:pr6}   are the main of them and they      
assert that, under very general conditions, a vector lattice either 
has a singleton d-basis of else this d-basis must be infinite. This
extends some of our work in [AK2, Section~6].

To make the reading of the article as much independent of [AK2] as       
possible we  collect in the next section  some necessary
definitions  and facts  about d-bases.  
Most of this preliminary material, as well as
some appropriate history regarding the subject, can be 
found in [AK2].

\bigskip
\section{Some preliminaries regarding d-bases} \label{prelim}

In terminology regarding vector lattices we follow [AB]. All vector     
lattices in this work are assumed to be Archimedean.
Whenever $B$ is a projection band in a vector lattice we denote by     
$P_B$ or $[B]$ the band-projection on $B$.

Recall [AK2, Definition~4.2] that a vector
lattice $X$  has a {\bf cofinal family of projection bands} if
for each non-zero band $B$  in $X$ there is a non-zero projection band    
$B_1\subseteq B$. Under a different name the same concept
was originally introduced in [LZ,~Definition~30.3].

\begin{definition} \label{d:d-basis}   %\marginnote{d:d-basis}
Let $X$ be a  vector lattice with a cofinal family
of projection bands. A collection of vectors 
$\{e_\gamma:\ \gamma \in \Gamma\}\in X$ 
is said to be {\bf d-independent} if for each projection band $B$ in $X$
the set $\{P_Be_\gamma: P_Be_\gamma\neq 0, \,  \gamma \in \Gamma\}$
is linearly independent, that is, the collection of all non-zero     
projections of the elements $e_\gamma$ on $B$ is linearly independent.
Any maximal $($by inclusion$)$ set of d-independent vectors is called a     
{\bf d-basis}.
\end{definition}

\medskip

A straightforward application of Zorn's lemma 
shows that in any vector lattice with a cofinal family of projection bands
(in particular, in any Dedekind complete vector lattice and in
any vector lattice with the projection property) there exists a
d-basis.

We will explain now what type of  
representation the  elements 
in $X$ have when a d-basis  $\{e_\gamma\}_{\gamma \in \Gamma}$          
is fixed. Namely, for each $x\in X$ there is a full collection
$\{X_i\}_{i\in I}$  of pairwise disjoint projection bands  
(depending on $x$) such that for each index $i$ the set 
$\Gamma_i= \{\gamma\in \Gamma: [X_i]e_\gamma\neq 0\}$ 
is finite and the element  $[X_i]x$  is a linear
combination of these linearly independent  projections
$[X_i]e_\gamma $ with $\gamma\in \Gamma_i$,
i.e., for some scalars   $\lambda_\gamma^{(i)}$ we have
$$ 
[X_i]x =\sum_ {\gamma\in\Gamma_i}\lambda_\gamma^{(i)}[X_i]e_\gamma. \eqno(1)
$$        
It is precisely the possibility of such a
representation that justifies our use of the term ``basis" here.
We would like to stress that there is a drastic difference between the
concepts  of a Hamel basis and a d-basis. For instance, the cardinality      
of  the latter can be much smaller. In extreme cases a d-basis may have
only one element.

\medskip
To be able to utilize  d-bases  more effectively we need to recall
the operation  of the  complete union that is defined in               
[V, Chapter~4] as follows: {\sl if $(x_i)$ is a collection of         
pairwise  disjoint elements in  a vector lattice $X$ such that            
there exist $\sup x_i^+ $
and
$\sup_i x_i^-$ in $X$, then the element $\sup_i x_i^+ -\sup_i x_i^-$
is called the {\bf complete union} and is denoted by
${\mathrel{\mathop {\bf S}\limits_i}}x_i$.} In particular, if each           
$x_i \in X_+$, then ${\mathrel{\mathop {\bf S}\limits_i}}x_i$ 
coincides with $\sup_i x_i$.

When combined with  operation 
$\mathrel{\mathop  {\bf S}}$,  
representation (1) described above gives the following.
Fix an arbitrary $x\in X$. As said before we can find  a full         
collection $\{X_i\}$ of pairwise disjoint projection bands in        
$X$ and  some scalars $\lambda^{(i)}_\gamma$  (all depending on $x$)
such that for each $i$ only a finite number of coefficients 
$\lambda^{(i)}_\gamma$ may be non-zero and for them
$[X_i]x =\sum_{\gamma}\lambda_\gamma^{(i)}[X_i]e_\gamma$.
Since  the bands $\{X_i\}$ are pairwise disjoint and full in $X$
we necessarily have $x={\mathrel{\mathop {\bf S}\limits_i}}[X_i]x$,
and hence the following ``global" representation holds:
$$ 
x={\mathrel{\mathop {\bf S}\limits_i}}\Sigma_{\gamma\in \Gamma}
\lambda^{(i)}_\gamma [X_i]e_\gamma. \eqno (\star) 
$$
We will refer to representation $(\star)$ as a {\bf d-expansion} of         
$x$ (with respect to the d-basis $\{e_\gamma\}$). 

Formally speaking, a d-expansion is not unique since we can always    
subdivide any projection band $X_i$ into the direct sum of two        
complementary projection bands  (assuming, of course, that the band $X_i$ is
not one-dimensional). Essentially, however,  any d-expansion is  unique in
the following sense. If 
$x={\mathrel{\mathop {\bf S}\limits_j}}\Sigma_\gamma
\lambda^{(j)}_\gamma [B_j]e_\gamma$ is another d-expansion with          
a generating collection of mutually disjoint bands $\{B_j\}_{j\in J}$, then
necessarily $\lambda^{(i)}_\gamma = \lambda^{(j)}_\gamma$ whenever 
$[X_i\cap B_j] e_\gamma \neq 0$.

If we do not assume that
$X$ has  a cofinal family of band-projections, then we cannot
guarantee any longer the existence of a sufficient quantity of
band-projections  $[B]e_\gamma$ and, as a result of it, we loose the
possibility           of having  a very useful d-expansion $(\star)$. We will
show in Section~4 how to retain an analogue of the d-expansion 
for an arbitrary vector lattice $X$.

\section{Characterization of band preserving projections}

Representation $(\star)$ was crucial 
in [AVK] for producing an example of a band preserving       
operator that was not a band-projection. Recently, 
in [AK2],   modifying slightly the same idea, 
we have  improved the previous example by
constructing   a band preserving {\bf projection}      
operator that is not a band-projection.
In view of the central role played   by  such projection 
operators in many situations it seems desirable 
to get a deeper understanding of their structure. And, as it turns out, 
the language of d-bases is adequate for obtaining a complete
description of all band preserving projection operators on Dedekind
complete vector lattices.

A vector sublattice $X_0$ of a vector lattice $X$            
is called {\bf component-wise closed} in $X$ if for each $u\in X_0$
the set ${\mathcal C}(u)$ of all components  of $u$ in $X$ is a subset 
of $X_0$.
Recall that ${\mathcal C}(u)=\{v\in X: |v|\land |u-v| =0\}$.

Observe that each component-wise closed vector sublattice of a
vector lattice with the principal projection property (resp. 
the  projection property)  also satisfies the             
principal projection property (resp. the  projection property).
In particular, if $X_0$ is a component-wise closed vector sublattice         
of a Dedekind complete vector lattice, then $X_0$ has a cofinal family
of band-projections.

It is proved in Proposition~4.9 in [AK2] that
{\sl if $T$ is a disjointness preserving operator from a
vector lattice $X$ into a vector lattice $Y$, then for each
ideal $Y_0$ in $Y$ its inverse image $X_0=T^{-1}(Y_0)$
is a vector sublattice of $X$ and, moreover, $X_0$ is 
component-wise closed in $X$.}

A simple but useful Lemma~7.3 in [AK2] asserts
that {\sl if $X$ is  a  vector lattice with 
a cofinal family of band-projections, then
a linear operator $T:X\to X$ 
is band preserving if and only if 
it commutes with band-projections, that is, $TP=PT$ for each  
band-projection $P$ on $X$.}

These two results will be used in our 
characterization  of band preserving projection
operators on laterally complete (in particular, universally complete)    
vector lattices. But first we  present a 
useful characterization of band-projections.

\begin{proposition} \label{p:band-pr}  %\marginnote{p:band-pr}

Let $X$ be a  vector lattice and $P$ be a band preserving
projection operator on $X$. The operator $P$ is a band-projection             
if and only if its kernel $P^{-1}(0)$ is a projection band.
\end{proposition}

\begin{proof}
Only the ``if" part is non-trivial.
Assume that $A=P^{-1}(0)$ is a projection band in $X$.
Consider the complimentary band $B=A^d$. Then $A\oplus B=X$.
Since $P$ is band preserving it follows that $P$ leaves $B$
invariant. But then the restriction of the operator  $P$ to $B$              
is a projection operator with  the trivial kernel, and so $P$ 
is the identity on $B$.
\end{proof}

Recall that a vector lattice  is {\bf laterally complete} 
if each collection of pairwise disjoint elements has a supremum.
By a well known theorem of Veksler and Geyler [VG,~Theorem~8]                
each laterally complete vector  lattice necessarily has  the projection
property. 
The most important example of laterally complete vector lattices
is provided by the class of universally complete vector lattices.

\begin{theorem}  \label{t:char.pr.oper}   %\marginnote{t:char.pr.oper}
Let $W$ be a laterally complete vector lattice. 
There is a one-to-one correspondence between
band preserving projection operators on
$W$ and  vector sublattices of
$W$ which are {\bf component-wise closed}  and {\bf laterally complete}.
This correspondence is given by the map $P\mapsto P^{-1}(0)$,
where $P$ is an arbitrary band preserving projection 
operator on $W$.
\end{theorem}

\noindent
{\bf Proof.} 
Let $P$ be a band preserving projection operator  on $W$. 
Since $P$ preserves dis\-jointness,  Proposition~4.9          
in [AK2] cited above  implies that the kernel $W_0 = P^{-1}(0)$ is a
component-wise closed  vector sublattice of $W$.
The assumptions that $W$ is laterally complete and that $P$ 
is  band preserving imply immediately that $W_0$ is also 
laterally complete.

Conversely, let $W_0$ be a component-wise closed and laterally complete
vector sublattice of $W$. 
Hence $W_0$  has the projection property.
Therefore, as noted after Definition~\ref{d:d-basis},    
there is a  d-basis $\{u_\alpha\}$ in $W_0$.  
{Using  that $W_0$ is  component-wise closed we can easily see that 
$\{u_\alpha\}$ remains d-independent in $W$. Hence,}
by  Zorn's lemma, we can 
extend $\{u_\alpha\}$ to a d-basis in $W$, that is,
to find  elements $\{v_\beta\}$ in $W$ such that 
$\{u_\alpha \}\cup \{v_\beta\}$            
is a d-basis in $W$. Now we are ready to define a necessary
operator $P$  on $X$.
Take an arbitrary $x\in W$ and consider its d-expansion
$(\star)$ with respect to the d-basis $\{u_\alpha \}\cup \{v_\beta \}$:
$$ 
x={\mathrel{\mathop {\bf S}\limits_i}}
\Big(
\Sigma_{\alpha}\lambda^{(i)}_{\alpha} [X_i]u_\alpha + 
\Sigma_{\beta} \lambda^{(i)}_{\beta} [X_i]v_\beta \Big). 
$$
The image $Px$ is defined by ``ignoring" the contribution of the
first part of the d-basis, that is,
$$ 
Px:={\mathrel{\mathop {\bf S}\limits_i}}
\Sigma_{\beta} \lambda^{(i)}_{\beta} [X_i]v_\beta. 
$$
Since $W$ is laterally complete and the elements 
\ $(\Sigma_{\beta}\lambda^{(i)}_{\beta} [X_i]v_\beta)_i$ \ 
are pairwise disjoint
their complete union (that is, the  value of $P$ at $x$) exists.
We omit a straightforward verification
that $P$ is a well defined band preserving projection operator on        
$W$ and that $\ker(P)=W_0$. 
It is worth pointing out that  
the assumption that $W_0$ is laterally complete
is essential  for
the validity of
the equality $\ker(P)=W_0$.

     {Note that $I-P$ is   also a band preserving
projection operator and so, by the first part of the theorem,
the kernel $ker(I-P)$ or, equivalently, the range of $P$                     
is a component-wise closed laterally complete vector sublattice 
of $W$. Hence, there are  d-bases in $ker(I-P)$.}

It remains to prove that if $P_1$ and $P_2$ are two band              
preserving projection operators such that $\ker (P_1)=\ker(P_2)$, then        
$P_1 = P_2$. This is a very special feature of band preserving  projections
that does not hold for general projections.
Let $\{u_\alpha \}$  be a d-basis in $\ker (P_1)$ 
and $\{v_\beta \}$ be a d-basis in $\ker (I-P_1)$.
Let us  verify that the system 
$\{u_\alpha \}\cup \{v_\beta\}$ is d-independent. 
Assume that for some band $B$ in $W$ we have
$$ 
\sum_i\lambda_i [B]u_{\alpha_i} + \sum_j\mu_j [B]v_{\beta_j} =0 \eqno (2) 
$$
Let us apply $P_1$ to this identity, keeping in mind that $P_1$
commutes with $[B]$ and that $P_1(u_\alpha) =0$ for each $\alpha$ 
and $P_1(v_\beta) =v_\beta$ for each $\beta$. We obtain that
$$ 
 \sum_j\mu_j [B]v_{\beta_j} =0, 
$$
implying that  $\mu_j=0$ whenever $[B]v_{\beta_j}\neq 0$
since $\{v_\beta\}$ are d-independent. Similarly,
applying to (2)  operator $I-P_1$ we will conclude
that $\lambda_i=0$ whenever $[B]u_{\alpha_i}\neq 0$.

Now, the
identity  $P_1x + (I-P_1)x =x $ implies  that
$\{u_\alpha \}\cup \{v_\beta\}$ is, in fact,  a d-basis in $W$. 

Take any  $z\in \ker (I-P_2)$. Since $\ker (P_1)=\ker (P_2)$ it follows
(similarly to our arguments above) that  the
system $\{u_\alpha, z \}$ is d-independent. Therefore $z$ allows a
d-decomposition with respect to $\{v_\beta \}$ and hence 
$z\in \ker (I-P_1)$. Thus we see that 
$\ker (I-P_2)\subseteq \ker
(I-P_1)$, and  similarly $\ker (I-P_1)\subseteq \ker (I-P_2)$.                
This establishes that $P_1=P_2$.
\qed

\medskip

Recall [AK2]  that a Dedekind complete vector lattice $X$ is
said to be {\bf principally universally complete} 
if each principal band in $X$ is 
universally complete.

\begin{corollary}   \label{c:puc}    %\marginnote{c:puc} 
Let $X$ be a principally universally complete vector lattice.
There is a one-to-one correspondence between band preserving projections on    
$X$ and vector sublattices $X_0$ of $X$ satisfying the following two
conditions:

{\rm 1)} $X_0$ is component-wise closed  and

{\rm 2)} For each  principal band $B$ in $X$ the intersection
 $B\cap X_0$  is laterally complete.
\end{corollary}

Now we are able to give a complete description of band preserving
projections on arbitrary Dedekind complete lattices.

\begin{theorem} \label{t:main.thm}  %\marginnote{t:main.thm}

Let $X$ be a Dedekind complete vector lattice and $P$ be a band preserving
projection operator on $X$.
Then $X=X_1 + X_2$, where the uniquely determined complimentary bands  
$X_1$ and $X_2$ satisfy the following properties.

{\rm 1)} $X_1$ is the maximal band such that the restriction of $P$ to $X_1$ 
is a regular operator, and therefore $P|X_1$ is a band-projection.

{\rm 2)} $X_2$ is principally universally complete, and hence the restriction  
$P|X_2$ is described by the previous corollary.

\end{theorem}
\begin{proof}
Let $X_1$ be the maximal band of regularity of $P$.
The existence of this band was established by de Pagter [P],
and has been  reproved in  Theorem 14.8 in [AK2]. The latter
theorem  asserts also that
the complimentary band $X_2=(X_1)^d$ is
principally  universally complete. An application of Corollary~\ref{c:puc} 
to the band $X_2$  finishes  the proof.
\end{proof}

Comparing Theorems~\ref{t:char.pr.oper} and~\ref{t:main.thm}                
one is led to ask whether the latter theorem can be generalized              
to vector lattices with the projection property. Somewhat unexpectedly
the answer to this question is negative. Let us explain why.
Theorem~\ref{t:main.thm}
implies  that {\it if a Dedekind complete vector lattice $X$ does             
not have a non-trivial principally universally complete band, then          
each band preserving projection operator  on $X$ is a band-projection}.
Therefore if we can produce an example of a 
band preserving projection operator $P$ on a
normed vector lattice $X$ such that
i) $X$  has the projection property and ii) $P$ is not a
band-projection, then this will establish that a generalization                 
in question cannot be true. (Keep in mind that no normed vector       
lattice  can have a non-trivial principally universally complete band.)    
An example like that can be  easily constructed directly or by modifying
slightly the example given in [AK1, Theorem~2] and [AK2, Theorem~13.1]. 
A less direct proof that such a generalization is not possible
is as follows.
If it
were true, then  on each vector lattice
without  non-trivial principally universally complete bands
each band preserving projection operator would be a band-projection.
In other words,  each  vector lattice $X$ like that  would have       
a determining family of band-projections [AK2, Definition~7.2].              
But then, by Theorem~8.5 in [AK2], each disjointness preserving
bijection $T$ from
$X$ onto an arbitrary vector lattice $Y$ with a cofinal family of
band-projections would have a disjointness preserving inverse.             
This  contradicts [AK1,~Theorem~2] and [AK2, Theorem~13.1].

Thus, we have shown that
Theorem~\ref{t:main.thm} cannot be generalized to vector
lattices  with the projection property. In other words, without the
$(r_u)$-completeness of the vector lattice  the projection              
property alone is not enough for the validity of Theorem~\ref{t:main.thm}.

\section{d-independence}

All previous work [AVK, AAK, AK2] on d-bases, in particular the definitions
of d-independence and of d-basis, as introduced in Definition~\ref{d:d-basis},
depends heavily on the availability of sufficiently many projection bands
in the corresponding vector lattices. This means that
these definitions cannot be easily  adopted to arbitrary vector lattices.
At the same time, the importance of the results 
devoted to d-bases and their numerous applications
suggest that it would be
desirable to extend these    results to arbitrary vector lattices. 
To do so one has to avoid using projection bands in the               
definition of d-independence. As we will see, this is possible
and  will be done in  this section.

\begin{definition}            \label{d:d-ind}       %\marginnote{d:d-ind}

A  system of elements
$\{x_\gamma\}_{\gamma\in \Gamma}$ in a vector lattice $X$ is 
called d-inde\-pendent if for every band $B$ in $X$, for every finite
set  of indices 
$\{\gamma_1,\ldots,\gamma_n\} \subseteq  \Gamma$, and every finite set of
non-zero scalars 
$\{c_1,\ldots, c_n\}$ the following implication holds:
$$
{\rm If} \ \ \sum\limits_{j=1}^nc_jx_{\gamma_j}\perp B,\ \ {\rm then} \ \ 
x_{\gamma_j}\bot B\ \ {\rm for \ each} \ j=1,2,\ldots, n.
$$
\end{definition}

\noindent
Replacing  each band $B$ above by its disjoint      
complement $E=B^d$, we obtain equivalently that a system
$\{x_\gamma\}_{\gamma\in
\Gamma}$ is d-independent if and only if for every band $E$ in
$X$,                       for every finite set  of indices 
$\{\gamma_1,\ldots,\gamma_n\} \subseteq  \Gamma$, and                      
every finite set of non-zero scalars 
$\{c_1,\ldots, c_n\}$ the following implication holds:
$$
{\rm If} \ \ \sum\limits_{j=1}^nc_jx_{\gamma_j}\in E,\ \ {\rm then} \ \ 
x_{\gamma_j}\in E\ \ {\rm for \ each} \ j=1,2,\ldots, n.
$$
It is easy to see that {\bf for any vector lattice $X$ with a cofinal
family of projection bands the above definition of
d-independence is equivalent to that given in 
Definition~\ref{d:d-basis}}. Accordingly, we retain the term.
To illustrate the difference between Definitions~\ref{d:d-basis} 
and~\ref{d:d-ind} consider the following example. 
Let $X=C[0,1]$. 
Since there is no non-trivial projection band  in $C[0,1]$, any               
two linearly independent functions $x_1,x_2 \in X$ are d-independent          
in the old sense. However, this two functions can easily be d-dependent       
in the sense of our new definition. Indeed, if there is a point  
$t_0\in [0,1]$ such that $x_1(t)=x_2(t)$  for all $t$ in a vicinity of    
$t_0$, then clearly  $x_1$ and $x_2$ are not d-independent in             
the sense of Definition~\ref{d:d-ind}.

More generally, a collection of functions $\{x_i\}_{i=1}^m$
in a vector lattice $X=C(K)$, where $K$ is a compact Hausdorff space, 
fails to be
d-independent if and only if there is a non-empty open subset $U$ of $K$
and scalars $\{c_i\}_{i=1}^m$, not all zero, such that
$\sum_{i=1}^m  c_ix_i(t) =0$ for each $t\in U$. This leads us to the following
result that will be used repeatedly.

\begin{lemma}  \label{l:polyn} %\marginnote{l:polyn}

For a function
$x\in C(K)$, where  $K$ is a compact Hausdorff space, the following statements
are equivalent:

{\rm 1)} For each $m\ge 2$ the collection $x,x^2,\ldots,x^m$ is
d-independent.

{\rm 2)} For some  $m\ge 2$ the collection $x,x^2,\ldots,x^m$ is
d-independent.

{\rm 3)} For each non-empty open $U\subseteq K$ the restriction of $x$ to $U$ 
is not a constant.

\end{lemma}

\begin{proof} The implications $1) \Rightarrow 2) \Rightarrow 3)$
are obvious. To prove $3) \Rightarrow 1) $ assume, contrary
to what we claim, that there is some $m\ge 2$ for which the           
collection $x,x^2,\ldots,x^m$ is not d-independent. Then
there is a non-empty open subset $V$ of $K$ and some scalars $c_i$ not all of
which are zero such that
$\sum_{i=1}^m  c_ix^i(t) =0$ for all $t$ from $V$.
Note, however, that $\sum_{i=1}^m  c_ix^i(t) =0$ is an algebraic equation
of degree not exceeding $m$. Therefore it cannot have more than $m$
solutions,  and this implies that the function $x$ must be a constant on a
non-empty open subset $U$ of $V$, a contradiction.
\end{proof}

\begin{corollary}  \label{c:polyn} %\marginnote{c:polyn}

For a function
$x\in C(K)$ and a non-empty open subset $U$ of   $K$ 
the following statements are equivalent:
\begin{itemize}
\item[{\rm 1)}] For each $m\ge 2$ the collection 
of the restricted to $U$ functions
$x|_U,x^2|_U,\ldots,x^m|_U$ is
d-independent in the space $C(U)$.

\item[{\rm 2)}] For some  $m\ge 2$ the collection 
of the restricted to $U$ functions
$x|_U,x^2|_U,\ldots,x^m|_U$ is
is d-independent in the space $C(U)$.

\item[{\rm 3)}] For each non-empty open $V\subseteq U$ the restriction of    
$x$ to $V$  is not a constant.
\end{itemize}
\end{corollary}

\begin{corollary}  \label{c:polyn2} %\marginnote{c:polyn2}

If a function $x\in C(K)$ satisfies 
Statement {\rm 3)} of the previous corollary
on a non-empty open subset $U$ of $K$, then
for each $m,k\ge 1$ and any scalars $\{\alpha_i\}_{i=1}^m $
and $\{\beta_i\}_{i=1}^k $ the functions 
$(\sum_{i=1}^m \alpha_i x^i)|_U$ and $(\sum_{i=1}^k \beta_i x^i)|_U $ are
d-dependent in the space $C(U)$ if and only if  $m=k$ and $\alpha_i=c\beta_i$
for some non-zero scalar $c$ and for each $i$.
\end{corollary}

The next statement is an immediate consequence of  Zorn's lemma.

\begin{proposition}  \label{p:d-max} %\marginnote{p:d-max}

Let $X$ be a vector lattice. There exists a maximal $($by inclusion$)$ 
d-independent system in $X$.
\end{proposition}

As surprisingly as it may sound, in comparison with
the case of  vector lattices with a cofinal family of
band-projections, it is  unclear yet whether or not
every  maximal d-independent system in a 
vector lattice $X$ is  sufficient to obtain an
analogue of  
d-expansion $(\star)$ 
for elements  in $X$. 
As we will show below in many cases it is sufficient. 
But in general we need to introduce a formal definition of a d-basis
in an arbitrary vector lattice.

\begin{definition}  \label{d:d-basis-new}   %\marginnote{d:d-basis-new}

We say that  a d-independent system $\{e_\gamma\}_{\gamma\in\Gamma}$
in a vector lattice $X$ is a {\bf d-basis} 
if for each $x\in X$ we can find a full system 
$\{X_i: i\in I\}$ of pairwise disjoint bands in $X$ 
such that for each $i \in I$
there is a finite number of indices $\gamma_1, \ldots,\gamma_n$
and non-zero scalars $c_1, \ldots, c_{n}$ $($all depending on $x$ and $i)$
such that  
$$
x-\sum\limits_{j=1}^{n} c_je_{\gamma_j}\perp X_i.  \eqno(\star\star)
$$
\end{definition}

Obviously each d-basis is a maximal d-independent system but
the validity of the converse statement is not known yet. 
If each band $X_i$ in the previous definition is a projection
band,  then  identity $(\star\star)$ can be rewritten  as 
$$
[X_i]x= \sum\limits_{j=1}^{n}c_i[X_i]e_{\gamma_j},
$$
which is nothing else but  d-expansion $(\star)$. 
That is, we see that in the case {\bf when $X$ has a cofinal
family of projection bands the new  definition of a d-basis           
coincides  with the old one}.

\bigskip
\noindent
{\bf Problem 1.}
Is every maximal d-indepen\-dent system  a d-basis?

\medskip

If the answer to this problem is negative, then
the following two open problems are of 
considerable interest.

\bigskip
\noindent
{\bf Problem 2.}
Describe the class of vector lattices  admitting a d-basis.

\bigskip
\noindent
{\bf Problem 3.}
Describe the class of vector lattices in which every
maximal d-indepen\-dent system is a d-basis.

\medskip

Below we will introduce a condition under which the answer to       
Problem~1 is affirmative. The condition is rather            
technical but, nonetheless, it holds for many important classes of       
vector  lattices.

\begin{definition}   \label{d: semi-c}

Let $x$ be a non-zero element in a vector lattice $X$.
We say that a non-zero element $b\in X$ is a {\bf semi-component} of $x$ 
if there is a full in $X$  system of pairwise disjoint bands  
$\{B_j: \ j\in J\}$   and a system of scalars  $\{c_j: \  j \in J \}$  
such that $b-c_j x\perp B_j$ for each $j\in J$.
\end{definition}

Formally speaking, $b=0$ can be also considered as a semi-component,         
but we are  excluding  this trivial case. Each non-zero
component of $x$  is obviously a semi-component as well.                    
This shows that the notion of semi-components generalizes                   
that of components. 

As we will see below, it is of special
importance when a given element has a semi-component in a given
band. 
So, consider  a  band $B$ in $X$ and let $x\notin B^d$.
If there is a non-zero component $b$ of $x$ that belongs to $B$,
then  $b$ is a semi-component of $x$ in $B$. In particular,
the band-projection $b=[B]x$, whenever it exists, is a semi-component of    
$x$ in $B$.

In general, it is easy to
verify                
that a non-zero element $b\in B$ is a semi-component of $x$
if there is a full in $B$  system of pairwise disjoint bands  
$\{B_j: \ j\in J\}$   and a system of scalars  $\{c_j: \  j \in J \}$  
such that $b-c_j x\perp B_j$ for each $j\in J$. 

Note also that the set of all semi-components of $x$ in $B$,              
together with the zero vector, is a vector sublattice in $X$.

\begin{example}  \label{e:semi-comp}
We describe here an important type  of vector lattices 
that have plenty of semi-components in each band. At the same time,
some of these vector lattices do not have a single non-trivial projection
band.
\end{example} 

Consider any compact space $K$ such that $X=C(K)$ has a dense
subspace of essentially constant functions. 
(Following {\rm [AK2]} we say  that a  function  $f\in C(K)$ is              
{\bf essentially constant} if for each non-empty open set
$G\subseteq K$ there exists a non-empty open subset $G_1\subseteq G$        
such that $f$ is  constant on $G_1$.)
For instance, each  metrizable compact space  $K$ or, more generally,            
each compact space with the countable chain condition,  satisfies this 
property (see [AB, Theorem~12.2], [HK, Theorem 0.1],  
[RR, Proposition, p.~130]).

Observe  next that  for an arbitrary non-zero band $B$ in $C(K)$
there is a  non-zero essentially  constant       
function $f$ in $B$. Indeed, fix an arbitrary 
function $g \in B$ such that $\|g\|=1$
and  ${\bf 0}\le g\le \bf 1$. 
By the hypothesis on $K$  there is an essentially             
constant function $g'\in C(K)$ such that $\|g-g'\|<\varepsilon$
for some small $\varepsilon>0$.
Consider $f=(g'-\varepsilon {\bf 1})^+$.  We omit a trivial
verification that $f$ 
is essentially constant, non-zero  and belongs to $B$.

To show that there are semi-components in any   band $B$,
take any $x_0\notin B^d$ and take an arbitrary non-zero             
essentially constant  function $f\in B$. The function  $fx_0$ 
belongs to $B$ and is a semi-component of $x_0$. 

If, additionally, the compact space $K$ is connected, 
then $X$ has no non-trivial projection band.
\hfill\qed

\begin{definition}   \label{d: new.cond}
We will say that a vector lattice $X$ satisfies 
condition 
$(*)$ if for every band $B$ in $X$ and every
$x \notin B^d$ there exists a semi-component of $x$ in $B$.

\end{definition}

All vector lattices with the projection property, or with the principal    
projection property, or just with a cofinal family of band-projections  
satisfy $(*)$. Moreover, if a vector lattice $X$ is such that {\it
for each band $B$ in $X$ and for each element       
$x\notin B^d$ there exists a non-zero {\bf component} 
of $x$ that belongs to $B$}, then $X$ also satisfies     
$(*)$. The class of  vector lattices with this latter property contains
properly the class of vector lattices  with a cofinal family of
band-projections. As we explained above, all vector lattices $C(K)$,       
where $K$ is a metrizable compact space or a compact space with the    
countable chain condition, also  satisfy  $(*)$.

\begin{proposition} \label{p:two.prob} %\marginnote{p:two.prob}

If a vector lattice $X$ satisfies $(*)$ then each
maximal $d$-indepen\-dent system in $X$ is a $d$-basis.
\end{proposition}

\begin{proof}
Let $\{x_i\}$ be a maximal d-independent system in $X$.
Assume, contrary to what we claim, that $\{x_i \}$ is not a d-basis. 
Then there exists an element $x\in X$ that cannot be 
d-expanded with respect to $\{x_i\}$ in the sense of
Definition~\ref{d:d-basis-new}. That is, we cannot find a full in $X$
collection of bands $X_i$ satisfying $(\star\star)$. 
This implies that there exists a
band
$B$ in
$X$ such that for any  non-trivial band  $B'\subset B$, 
for any finite set of indices $\{i_1,\ldots,i_n\}$, and for any
scalars $\lambda_1, \ldots, \lambda_n$ 
the element
$x-\sum\limits_{k=1}^n \lambda_k x_{i_k}$ is not disjoint to $B'$.
Let us consider this band $B$.
Since 
$X$ satisfies $(*)$, we can find a non-zero $b\in B$, a full
in $B$ system of pairwise disjoint bands $\{B_j\}$, and scalars $c_j$     
such that for each $j$ we have $b-c_j x \perp B_j$. It is
easy to verify now that the system $\{b, \ x_i \}$ is d-independent,     
a contradiction to the maximality of $\{x_i\}$.
\end{proof}

An important class of vector lattices satisfying condition $(*)$
is described next. Recall that for each $x$ in a vector lattice $X$
there exists a compact Hausdorff space $K_x$ such that
the principal ideal 
$X(x)=\{x'\in X: |x'|\le \lambda |x|, \ \lambda \in\R\}$ 
is order isomorphic to an  order dense vector sublattice of 
$C(K_x)$. We say that $X(x)$ is represented in $C(K_x)$.
This representation is unique
(up to a homeomorphism of  $K_x$) if we require            
that the element $x$ be mapped to the constant one function $\bf 1$.

\begin{theorem}  \label{t:r.impl.basis} 
Each $(r_u)$-complete vector lattice $X$ 
with the  countable sup property satisfies  condition $(*)$. 
In particular, every maximal d-independent system in $X$ 
is a d-basis.
\end{theorem}

\begin{proof}
To verify that  $X$ satisfies condition $(*)$
take any band $B$ in $X$ and any element $x\in X_+$
that is not disjoint to $B$.
Consider the principal ideal $X(x)$ generated by $x$.
Then $X(x)=C(K)$ for some compact Hausdorff space $K$
and $B\cap C(K)$ is a band in $C(K)$.
Since $X$ satisfies the  countable sup property,
the compact space $K$  satisfies the countable chain condition
and consequently,  as said earlier, the collection of essentially 
constant functions is dense in $C(K)$ in view of [HK,RR]. 
As shown in Example~\ref{e:semi-comp}
there is a  non-zero essentially  constant       
function $f$ in $B\cap C(K)$. It remains to note that
$f$ is a semi-component of $x$ not only in $C(K)$ but also in $X$,          
and we are done.
\end{proof}

\medskip

Now we are going to discuss some questions related to the cardinality of   
maximal d-independent systems and d-bases. More precisely we are     
interested in the following questions.

\bigskip
\noindent
{\bf Problem 4.}   
Let $X$ be a vector lattice. Do
all maximal d-independent systems in $X$ have the same cardinality?

\bigskip
\noindent
{\bf Problem 5.}
Let $X$ be a vector lattice admitting a d-basis. 
Do all  d-bases in $X$ have the same cardinality?

\medskip
Regardless of the answers to the previous problems, it will be also of  
interest to relate the cardinality of a maximal d-independent system      
(resp. of a d-basis) to some other cardinal characteristics of the vector
lattice $X$, for example, to the disjointness type $t(X)$ introduced in [AV].

\bigskip

If a vector lattice $X$ does not have a weak unit, then obviously any       
d-independent maximal system in $X$ is infinite. 
The  vector lattice $c_{00}$ of eventually zero sequences               
provides an example of  a  discrete Dedekind complete vector             
lattice without a weak unit. Every d-basis in $c_{00}$ is countable.

Discrete vector lattices with a weak unit provide the simplest example      
possible when there is a singleton d-basis. 
Discreteness, however, is not a decisive factor here.
Each essentially one-dimensional vector lattice [AK2] with                      
a weak unit has a singleton d-basis.
Moreover, Gutman [G] constructed an
example  of an extremally disconnected compact Hausdorff space $K$   
without isolated points  such that every continuous function on $K$ is
essentially constant. In other words, Gutman's space $C(K)$ is atomless
and, nevertheless,  has a singleton d-basis.   
Another example  of an atomless (but not Dedekind
complete) vector   lattice with a singleton    d-basis is the space 
$C(\beta \N\setminus \N)$.  Similar examples can be found in
[HK,RR].                   
It is an interesting  open problem  {\it to describe 
compact Hausdorff spaces $K$ for which the vector lattice  $C(K)$ has 
a singleton d-basis}. 
\medskip

If we do not impose any conditions on a vector lattice, then we          
can easily find a vector lattice $X$ with a d-basis of a pre-assigned 
finite cardinality. Perhaps, the simplest example of a vector lattice      
with a d-basis of cardinality $n,\ 1\leq n<\infty $, is  provided 
by the vector lattice $X$ of piecewise polynomials of degree
not exceeding $n$, that is, a continuous on [0,1] function $x$           
belongs to $X$ if and only if  there is a partition
$t_0= 0 < t_1< \ldots < t_{m-1} <1 = t_{m}$  of [0,1] and
polynomials $p_1, \ldots, p_{m}$ of degree no more than $n$ 
such that $x\equiv p_j$ on $[t_{j-1},t_j]$. This space was   
considered in [AK2].

If  we do not restrict the degrees of the  polynomials $p_j$ 
above,  
then we obtain an example of a vector lattice with a countable
d-basis. This vector lattice is a {\it subalgebra} of $C[0,1]$ and, as        
we shall demonstrate below, this additional algebraic structure 
is the actual reason  of why any d-basis  is infinite. 
\medskip

Two major results regarding cardinality were proved in [AK2].
Namely, in Theorem~6.8  we proved that a Dedekind complete vector         
lattice $X$ either has a singleton d-basis or else any d-basis in $X$ 
is infinite. After that          
we showed that for two large and  important classes of vector 
lattices any d-basis is, in actuality, uncountable. 
One of this classes consists of all non-zero ideals in the            
space $L_0[0,1]$ (Theorem~6.10 in [AK2]) and the other class 
consists of all non-zero ideals in the universal completion of the              
$C(K)$ space, where  $K$ is  any compact metric space 
without isolated points (Theorem~6.9 in [AK2]).

Below we will extend the first of these results to a much
more general class of vector lattices than Dedekind complete.
As for Theorems~6.9 and~6.10 in [AK2], it is presently unclear
whether or not each infinite d-basis in an atomless
Dedekind complete vector lattice must be  uncountable.

\begin{definition}
Let us say that a vector lattice $X$ is a {\bf local algebra} if
for each element $x\in X$ the principal ideal $X(x)$ can be
represented as a {\bf subalgebra} and a vector sublattice
of the corresponding space $C(K_x)$.
\end{definition}

There are many non-trivial examples of local algebras, for instance 
all $(r_u)$-complete vector lattices (in particular, all Dedekind        
complete vector lattices) are local algebras.  Furthermore,
a solution of the next problem can help find more such examples.

\bigskip
\noindent
{\bf Problem 6.}
Is it possible to describe  vector lattices that are 
local algebras in terms not involving representations?

\medskip

\begin{theorem}         \label{p:pr5}
Assume that a vector lattice $X$  is a local algebra with
no singleton d-basis. Then every d-basis in $X$, whenever it  exists, is
infinite.        
\end{theorem}

\begin{proof} 
Suppose that $e_1,\ldots,e_n$ is a d-basis in $X$, where $n\ge 2$.
There are at least two elements $e_i,\; e_j$, $1\leq i<j\leq n$,           
which are    not disjoint. Indeed, otherwise the element 
$e_1+ e_2+ ...+ e_n$ would be a singleton d-basis in $X$. 
Let $e=\sum\limits_{k=1}^n|e_i|$ and
let $\pi$ be an order isomorphism of $X(e)$ onto a  dense   
subalgebra of some $C(K)$ space such that $\pi e=\bf 1$. 
Since $\pi e_i$ and $\pi e_j$ are d-independent, there exists
a non-empty open subset $U\subset K$ such that at least one of the     
functions $\pi e_i,\; \pi e_j$ is not constant on each non-empty            
open subset $V\subset U$. Assume for  definiteness that the function        
$\pi e_i$ is such. Then for any positive integer $m$ the functions 
$\pi e_i,\ldots, (\pi e_i)^m$ belong to $\pi(X(e))$ and  by  
Corollary~\ref{c:polyn} they are d-independent on  $U$. 
To conclude the proof we will verify that this 
contradicts the              
fact that the system 
$\{\pi e_1,\ldots, \pi e_n\}$ 
is a d-basis in $\pi(X(e))$. Indeed, the latter implies that
for each non-empty open set $G\subseteq K$ and each $k=1,\ldots,m$
we can find a non-empty open subset of $G$ on which
the element $(\pi e_i)^k$ is a linear combination of the elements
of our d-basis. Consequently we can find a non-empty open subset $V$
of $K$
on which each of the  functions $\pi e_i,\ldots,(\pi e_i)^m$
is a linear  combination of the functions of the d-basis
$\{\pi e_1,\ldots, \pi e_n\}$. But this is impossible
if $m>n$ since, as said above, the functions 
$\pi e_i,\ldots,(\pi e_i)^m$
are linearly independent on $V$. 
\end{proof} 

\begin{corollary}         \label{c:pr5}
Let $X$ be an $(r_u)$-complete vector lattice.
Then either $X$ has a singleton d-basis or
every d-basis in $X$,  whenever it exists, is infinite.        
\end{corollary}

We do not know if  an analogue of the previous theorem remains         
true for   d-maximal systems instead of  d-bases. 
Under an additional condition we can prove this.

\begin{theorem}  \label{p:pr6} 
Assume that a vector lattice $X$  is a local algebra with
no singleton maximal d-independent system. If, additionally,   
$X$ satisfies the  countable sup property, i.e., every  family of
pairwise disjoint   non-zero elements in $X$ is at most countable,
then every maximal d-independent system in $X$ is infinite.
\end{theorem}

\begin{proof}
Suppose to the contrary that there is a finite maximal d-independent
system $\{x_1,\ldots,x_n\}$ with $n\geq 2$. Then
after representing $X(|x_1|+\ldots+|x_n|)$ on a compact space $K$,             
exactly as in the  proof of the previous theorem, we can find
at least one element in our system $\{x_1,\ldots, x_n\}$, let it be
$x_1$ for definiteness, and a non-empty open subset 
$U$ of $K$ such that for any open
$V\subset U$ and for any positive integer $m$ the functions 
$(\pi x_1)^j,\ 1\leq j\leq m$, are linearly independent on $V$. 
Since $X$ satisfies the countable sup property, the compact space $K$
satisfies the countable chain condition.

Fix an $m>n+1$. For each $m$-tuple
of scalars $\bar\alpha=(\alpha_1,\ldots,\alpha_m)$ consider the function
$$
y_{\bar\alpha}= \sum_{i=1}^m \alpha_i (\pi x_1)^i.
$$
Observe that if $\bar \alpha'=(\alpha_1',\ldots,\alpha_m')$
is another m-tuple. 
then by Corollary~\ref{c:polyn2} the functions
$y_{\bar\alpha}$ and $ y_{\bar \alpha'}$ are d-dependent on $U$
iff $\bar\alpha \equiv \bar c\alpha'$ for some scalar $c\neq 0$.

Let us fix an arbitrary  uncountable collection of points 
$\Lambda =\{\bar\alpha\} \subset \R^m$ 
such that any $n+1$ pairwise distinct points 
$\bar \alpha_1,\ldots, \bar \alpha_{n+1}$ in $\Lambda$ are linearly   
independent. Then, by Corollary~\ref{c:polyn2}, 
the functions $y_{\bar\alpha_1},\ldots, y_{\bar \alpha_{n+1}}$ 
are d-independent on $U$.

Fix any $\bar \alpha \in \Lambda$ 
and consider
the function $y_{\bar \alpha}$ introduced above. Since $\{x_1,\ldots,x_n\}$
is a maximal d-independent system, the function $y_{\bar\alpha}$
cannot be d-independent of this system on $U$. Therefore,
there is a non-empty open subset $U_{\bar\alpha} \subset U$
such that $y_{\bar\alpha}$  coincides on  
$U_{\bar\alpha}$ with a linear combination of the functions
$\pi x_1,\ldots,\pi x_n$. 

We claim next that whenever we take arbitrary $n+1$ mutually distinct
points $\bar \alpha_1,\ldots, \bar \alpha_{n+1} \in \Lambda$
we necessarily have that 
$$
\bigcap_{j=1}^{n+1} U_{\bar \alpha_j}=\emptyset  \eqno(3).
$$
Indeed, if an open set
$V=\cap_j U_{\bar \alpha_j}$ is not empty, then on $V$
each of the functions $y_{\bar\alpha_j}$ is a linear combination of
the functions from $\pi x_1,\ldots,\pi x_n$. On the other hand, as       
said    above, these functions $y_{\bar\alpha_j}$ are linearly       
independent on $V$, a contradiction.
 
We conclude the proof by showing that the existence 
of an uncountable family 
$\{U_{\bar\alpha} : \bar\alpha\in \Lambda\}$ 
of non-empty open sets satisfying (3), 
contradicts  the countable chain condition in $K$. 
We will use induction on $n$.

Assume that  the above statement is true for some $n$,  
and let us prove it for $n+1$. Consider all sets
of the form 
$$
V_\beta:=U_{\bar\alpha_1}\cap\ldots \cap U_{\bar\alpha_{n+1}}.
$$
where $\beta=(\bar\alpha_1, \ldots,  \bar\alpha_{n+1})$
and $\alpha_1,\ldots, \alpha_{n+1}$ are pairwise distinct points in 
$\Lambda$. Observe that by (3) the sets $V_\beta$ 
are pairwise disjoint. 
If each $V_\beta=\emptyset$, then we are done in view of the
induction hypothesis. So some sets $V_\beta$ are non-empty,
and, since they are pairwise disjoint, 
there may be at most countably many of such sets.
Let $V_{\beta_1},\ldots, V_{\beta_k}, \ldots$ be all these non-empty sets. 
Since each 
$\beta_k=(\bar\alpha^{(k)}_1, \ldots,  \bar\alpha^{(k)}_{n+1})$
depends on a finite number of indices in $\Lambda$, there are         
uncountably many
$\alpha$-s that have not  been used. Consider all the sets
$V_\beta$ for which  $\beta$ depends on at least one of these unused 
$\alpha$-s. Again by the induction hypothesis, all
these $V_\beta$ cannot be empty. Take an arbitrary such $V_\beta$ that is not
empty. However, this $V_\beta$ must be disjoint from each $V_{\beta_k}$, and so
we have a contradiction.
\end{proof}

\vfill
\eject

\bigskip
\centerline {\bf References}
\bigskip

[AK1] Y. A. Abramovich and A. K. Kitover, A solution to a problem on
         invertible disjointness preserving operators,
         {\it Proc. Amer. Math. Soc.\/ \bf 126}$\,$(1998), 
         1501--1505.

[AK2] Y. A. Abramovich and A. K. Kitover,  {\it Inverses of Disjointness 
      Preserving Operators}, {Memoirs of the Amer. Math. Soc.}, 
      forthcoming.

[AV]  Y. A. Abramovich and A. I. Veksler, Exploring partially ordered     
      spaces by means of transfinite sequences, 
      {\it Optimizacija, \/  No. \bf 12}$\,$(1973), 8--17.

[AVK] Y. A. Abramovich, A. I. Veksler, and A. V. Koldunov,
          Operators preserving disjointness, {\it Dokl. Akad. Nauk USSR \/
          \bf 248}$\,$(1979), 1033--1036.

[AD] C. D. Aliprantis and O. Burkinshaw, {\it Positive Operators\/}, 
        Academic Press, New York \& London, 1985.

[G] A. Gutman, Locally one-dimensional $K$-spaces and
       $\sigma$-distributive Boolean algebras, {\it Sibirian Adv.
       Math. \/ \bf 5}$\,$(1995), 99--121.

[HK] J. Hart and K. Kunen, Locally constant functions, {\it Fund.      
         Math.\/ \bf 150}$\,$(1996), 67--96.

[LZ] W. A. J. Luxemburg and A. C. Zaanen, {\it Riesz Spaces I\/},
           North-Holland, Amsterdam, 1971.

[P] B. de Pagter, A note on disjointness preserving operators,
           {\it Proc. Amer. Math. Soc.\/ \bf 90}$\,$(1984), 543--549.

[RR] M. E. Rudin and W. Rudin, Continuous functions that are locally constant
     on dense sets, {\it J. Funct. Anal. \/\bf 133}$\,$(1995), 120--137.

[VG] A. I. Veksler and V. G. Geyler, Order and disjoint completeness of
      linear partially ordered spaces,  {\it Siberian Math. J. \bf 13}
       $\,$(1972),  30--35.

[V] B. Z. Vulikh, {\it Introduction to the theory of partially ordered
          spaces,\/ } Wolters-Noordhoff Sci. Publication, 
           Groningen, 1967.

\vskip 2cm

Y. A.  Abramovich                    \hskip 5.5cm       A. K. Kitover

Department of Mathematical Sciences \hskip 1.92cm    Department of Mathematics

IUPUI, Indianapolis, IN 46202    \hskip 3.3cm     CCP, Philadelphia, PA 19130

USA                               \hskip 8cm    USA

yabramovich@math.iupui.edu      \hskip 3.5cm   akitover@ccp.cc.pa.us

\end{document}